\documentclass[12pt,a4paper]{amsart}
\usepackage{amsmath,amsfonts,amssymb}

\newtheorem{theorem}{Theorem}[section]
\newtheorem{proposition}[theorem]{Proposition}
\newtheorem{lemma}[theorem]{Lemma}
\newtheorem{corollary}[theorem]{Corollary}

\newcommand{\set}[2]{\ensuremath{\{ #1 \>|\> #2 \}}}

\begin{document}

\title{A Converse to the Second Whitehead Lemma}
\author{Pasha Zusmanovich}
\address{}
\email{justpasha@gmail.com}
\date{first written April 19, 2007; last minor revision December 24, 2013.}
\thanks{\textsf{arXiv:0704.3864}; J. Lie Theory \textbf{18} (2008), 295--299}

\begin{abstract}
We show that finite-dimensional Lie algebras over a field of characteristic zero
such that the second cohomology group in every finite-dimensional module vanishes,
are, essentially, semisimple.
\end{abstract}

\maketitle

\section*{Introduction}

The classical First and Second Whitehead Lemmata state that the first, respectively second,
cohomology group of a finite-dimensional semisimple Lie algebra with coefficients in any
finite-dimensional module vanishes. It is natural to ask whether the converse is true.
For the first cohomology this is well-known: 

\begin{theorem}[A converse to the First Whitehead Lemma]\label{first}
Any fi\-nite-\-di\-men\-si\-o\-nal Lie algebra over the field of characteristic zero
such that its first cohomology with coefficients in any finite-dimensional module 
vanishes, is semisimple.
\end{theorem}

\begin{proof} (see, e.g., \cite[Chapter 1, Theorem 3.3]{encycl-3}).
Due to cohomological interpretation of modules extension, the vanishing
of cohomology with coefficients in any finite-dimensional module is equivalent 
to the complete reducibility of any finite-dimensional module. 
Assuming that a Lie algebra $L$ with such condition has a nonzero abelian ideal $I$ 
and considering the adjoint representation of $L$, we get that $I$ is a
direct summand of $L$. But then $I$ also should satisfy the same condition, 
a contradiction. 
\end{proof}

What about the second cohomology?

Somewhat surprisingly, this question seemingly was not addressed in the literature.
The aim of the present elementary note is to show that, essentially, the converse is true too -- 
finite dimensional Lie algebras over a field of characteristic zero such that
their second cohomology with coefficients in any finite-dimensional module vanishes,
are very close to semisimple ones. In fact, the proof is only a bit more
complicated then the proof of Theorem \ref{first} and readily follows from 
results already established in the literature.

One may look at this question from a somewhat different angle. Recently, an interesting
class of so-called strongly rigid Lie algebras was investigated in \cite{petit} -- 
namely, Lie algebras
whose universal enveloping algebra is rigid. It turns out that for such algebras,
the second cohomology in the trivial module and in the adjoint module vanishes. 
Considering these two modules as the most ``natural'' ones, one may wonder for which 
Lie algebras the stronger condition -- vanishing of the second cohomology with 
coefficient in any finite-dimensional module -- would hold.

What happens in the positive characteristic? As was shown independently by 
Dzhumadildaev \cite{dzhu-sbornik} and Farnsteiner -- Strade \cite{fs}, for any finite-dimensional
Lie algebra over a field of positive characteristic and any degree
not greater than the dimension of the algebra, there is a module with non-vanishing cohomology in that degree.
(In fact, the low degree cases interesting for us here were settled even 
earlier -- for the first degree cohomology by Jacobson \cite[Chap. VI, \S 3, Theorem 2]{jacobson} 
and for the second degree cohomology -- again by Dzhumadildaev \cite{dzhu-umn}).
So the answer is trivially void in this case -- one of the rare
cases when the answer (but not the proof!) in the positive characteristic turns out to 
be simpler. The reason for this is
the possibility to construct analogues of induced modules with desired cohomological properties
by means of various truncated finite-dimensional versions of the universal enveloping algebra.
In characteristic zero, such modules would be infinite-dimensional.

\begin{theorem}[A converse to the Second Whitehead Lemma]\label{second}
A fi\-nite-\-di\-men\-si\-o\-nal Lie algebra over a field of characteristic zero such that
its second cohomology with coefficients in any finite-dimensional module vanishes, 
is one of the following:
\begin{enumerate}
\item an one-dimensional algebra;
\item a semisimple algebra;
\item the direct sum of a semisimple algebra and an one-dimensional algebra.
\end{enumerate}
\end{theorem}

In the first section of this note we accumulate the needed results from the literature, and
in the second section we provide the proof. 

Our notations are standard. If $L$ is a Lie algebra, $Rad(L)$ denotes the solvable radical
of $L$.
If $V$ is an $L$-module, $V^L = \set{v \in V}{xv = 0 \text{ for any }x \in L}$ is a submodule
of $L$-invariants.
Throughout the note, the ground field $K$ assumed to be of characteristic zero
and all algebras and modules assumed to be finite-dimensional. When being considered
as a module over a Lie algebra, $K$ understood as a trivial module.

\section{Needed results}

\begin{proposition}[Dixmier]\label{codim1}
Let $L$ be a Lie algebra, $V$ be an $L$-module, and $I$ be an ideal of $L$ of 
codimension 1. Then for any $n\in \mathbb N$ and $x\in L \backslash I$,
$H^n(L,V) \simeq H^n(I,V)^x \oplus H^{n-1}(I,V)^x$.
\end{proposition}

\begin{proof}
This was proved in an equivalent form of a certain long exact sequence in 
\cite[Proposition 1]{dixmier} 
and stated without proof (in a more general situation when $I$ is a subalgebra of codimension 1) 
in \cite[Proposition 4]{dzhu-func}. 
In fact, this is an easy consequence of the Hochschild--Serre spectral sequence.

Let $L = I \oplus Kx$ for some $x\in L$, $x\notin I$ and consider the Hochschild--Serre spectral 
sequence abutting to $H^*(L,V)$ relative to the ideal $I$. As
$E_2^{pq} = H^p (L/I, H^q(I,V)) \linebreak[0]\simeq H^p(Kx, H^q(I,V))$, it is nonzero only for $p=0,1$,
in which cases $E_2^{0q} = H^q(I,V)^x$ and $E_2^{1q} \simeq H^q(I,V) / x H^q(I,V) \simeq H^q(I,V)^x$.
Hence the spectral sequence stabilizes at $E_2$, 
$H^n(L,V) \simeq E_{\infty}^{0n} \oplus E_{\infty}^{1,n-1} = E_2^{0n} \oplus E_2^{1,n-1}$,
and the result follows.
\end{proof}

\begin{proposition}[Dixmier]\label{nilpotent}
If $L$ is a nilpotent Lie algebra of dimension $>1$, then $H^2(L,K) \ne 0$.
\end{proposition}

\begin{proof}
In \cite[Th\'eor\`eme 2]{dixmier}, a much more general result is stated:
$dim\,H^n(L,V) \ge 2$ for any $0 < n \le dim\,L$ and any finite-dimensional $L$-module 
$V$ containing $K$.
It is proved by repetitive applications of Proposition \ref{codim1}.
\end{proof}

And last, we will need the following simple result which implicitly contained already in the foundational paper
\cite{hs} and which is explicitly proved, for example, in \cite[Lemma 1]{rozen}:

\begin{proposition}[Hochschild -- Serre]\label{e2}
Let $L$ be a Lie algebra represented as the semidirect sum $L = P \oplus I$ of a subalgebra $P$ 
and an ideal $I$, and $V$ be an $L$-module. 
Then the Hochschild--Serre spectral sequence abutting to $H^*(L,V)$ relative to the ideal $I$, 
stabilizes at the $E_2$ term.
\end{proposition}

\section{Proof of Theorem \ref{second}}

Finite-dimensional Lie algebras with the property that their second 
cohomology with coefficients in any finite-dimensional module vanishes, will be
called \textit{2-trivial}.

\begin{lemma}\label{semidirect}
Let $L$ be a 2-trivial Lie algebra represented as the semidirect sum $L = S \oplus I$ of a subalgebra 
$S$ and an ideal $I$. Then:
\begin{enumerate}
\item $S$ is 2-trivial;
\item $(H^2(I,K) \otimes V)^S = 0$ for any $S$-module $V$.
\end{enumerate}
\end{lemma}

\begin{proof}
Let $V$ be an $L$-module and consider the Hochschild--Serre spectral sequence abutting to $H^*(L,V)$
relative to $I$. By Proposition \ref{e2}, it stabilizes at $E_2$, hence all $E_2$ terms vanish.
We have $E_2^{20} = H^2(L/I, H^0(I,V)) \simeq H^2(S, V^I)$ and 
        $E_2^{02} = H^0(L/I, H^2(I,V)) \simeq H^2(I,V)^S$.

Choose $V$ as follows: let $V$ be an arbitrary $S$-module, and $I$ acts on $V$ trivially. Then 
$E_2^{20} \simeq H^2(S,V)$,
$E_2^{02} \simeq (H^2(I,K)\otimes V)^S$, and they vanish for any $S$-module $V$, 
what proves (i) and (ii) respectively.
\end{proof}

In particular case when the semidirect sum reduces to the direct sum, (i) shows 
that the direct summand of a 2-trivial Lie algebra is 2-trivial.
The last assertion may be proved also directly by using the K\"unneth formula, 
without appealing to Proposition \ref{e2}.

\begin{lemma}\label{lemma-codim1}
An ideal of codimension 1 in a 2-trivial Lie algebra is a direct summand.
\end{lemma}

\begin{proof}
Let $L$ be a 2-trivial Lie algebra and $I$ be an ideal of $L$ of codimension 1. Write 
$L = I \oplus Kx$ for some $x\in L$. Evidently $adx$ is an $x$-invariant 1-cocycle in $Z^1(I,I)$. 
As by Proposition \ref{codim1}, $H^1(I,I)^x$ embeds into $H^2(L,I) = 0$, this cocycle is a coboundary,
i.e. there is $z \in I$ such that $[y,x] = [y,z]$ for any $y\in I$. Replacing $x$ by $x^\prime = x-z$,
we get a direct sum decomposition $L = I \oplus Kx^\prime$, $[I,x^\prime] = 0$.
\end{proof}

\begin{corollary}\label{2dim}
A $2$-dimensional Lie algebra is not $2$-trivial.
\end{corollary}

\begin{proof}
By Lemma \ref{lemma-codim1}, a $2$-trivial $2$-dimensional Lie algebra is abelian. 
But for any abelian Lie algebra $L$ and any $n\in \mathbb N$, $H^n(L,K) = C^n(L,K)$, 
and hence any abelian Lie algebra of dimension $>1$ is not $2$-trivial.
\end{proof}

\begin{lemma}\label{h2}
Let $L$ be a 2-trivial Lie algebra. Then $H^2(Rad(L), K) = 0$.
\end{lemma}

\begin{proof}
The assertion is obvious in the case when $L = Rad(L)$ is solvable, so suppose $L$ is not solvable.
Let $L = S \oplus Rad(L)$ be a Levi--Malcev decomposition, 
where $S$ is a semisimple Malcev subalgebra.

By Lemma \ref{semidirect}(ii), $(H^2(Rad(L),K)\otimes V)^S = 0$ for any $S$-module $V$. 
Assume $H^2(Rad(L),K) \ne 0$ and take $V = H^2(Rad(L),K)^*$.
There is a canonical surjection of $S$-modules $H^2(Rad(L),K) \otimes H^2(Rad(L),K)^* \to K$. 
Since any extension of $S$-modules splits, $K$ contained in 
$H^2(Rad(L),K) \otimes H^2(Rad(L),K)^*$,
hence $(H^2(Rad(L),K) \otimes H^2(Rad(L),K)^*)^S \ne 0$, a contradiction.
\end{proof}

Now we are ready to prove Theorem \ref{second}.

Let $L$ be a 2-trivial Lie algebra. We shall reason by induction on the dimension of $L$.

If $[L,L] = L$, then $Rad(L)$ is nilpotent 
(in fact, this is the consequence of the Levi Theorem which in turn is the 
consequence of the Second Whitehead Lemma; see 
\cite{jacobson}, Chap. III, \S 9, Corollary 2). 
Note that $Rad(L)$ cannot be one-dimensional, as then $Rad(L)$ is one-dimensional representation 
of a semisimple Lie algebra, and hence is a trivial representation, what
contradicts the condition $[L,L] = L$.
If $dim\,Rad(L) > 1$, then by Proposition \ref{nilpotent},
$H^2(Rad(L),K) \ne 0$, while by Lemma \ref{h2}, $H^2(Rad(L),K) = 0$, a contradiction.
Hence $Rad(L) = 0$, i.e. $L$ is semisimple.

Let $[L,L] \ne L$. As any subspace of $L$ containing $[L,L]$ is an ideal of $L$,
$L$ contains an ideal $I$ of codimension 1. By Lemma \ref{lemma-codim1},
$L$ is the direct sum of $I$ and an one-dimensional algebra, 
and by Lemma \ref{semidirect}(i), $I$ is 2-trivial. 
By induction assumption, $I$ is either one-dimensional, or semisimple, or the direct sum
of a semisimple algebra and an one-dimensional algebra. 
In the first case $L$ is abelian 2-dimensional, a contradiction (Corollary \ref{2dim}).
In the third case $L$ is the direct sum of a semisimple algebra and a 2-dimensional 
abelian algebra.
By Lemma \ref{semidirect}(i), both direct summands are 2-trivial, the same contradiction again.
Hence the only possible case is when $L$ is the direct sum of a semisimple algebra and 
an one-dimensional algebra, what concludes the proof.

\section*{Acknowledgements}

I am grateful to Dietrich Burde who asked the question which prompted me to write 
this note, and to anonymous referee whose suggestions improved the presentation.

\end{document}

\title{Erratum to \\ ``A converse to the Second Whitehead Lemma''}
\author{Pasha Zusmanovich}
\address{
Institute of Mathematics and Statistics, University of S\~{a}o Paulo,
Brazil
}
\email{pasha@ime.usp.br} 
\date{April 27, 2014}

\maketitle

In \cite{z1}, a result from \cite{rozen} is used, claiming that the 
Hochschild--Serre spectral sequence abutting to the cohomology of the 
semidirect sum $L = S + I$ of Lie algebras $S$ and $I$  
($S$ is a subalgebra, $I$ is an ideal, i.e. $S$ acts on $I$) with coefficients
in an arbitrary $L$-module $V$, with respect to the ideal $I$, stabilizes at the 
$E_2$ page (\cite[Lemma 1]{rozen} and \cite[Proposition 1.3]{z1}). 
In the whole generality, this statement is false. In fact, as shown in 
\cite{barnes}, this spectral sequence can stabilize at arbitrarily large step.

This statement is true, however, in many special cases: for example,
if the base field is of characteristic zero, $S$ is a 
finite-dimensional semisimple Lie algebra, and $V$ is finite-dimensional 
(\cite[Theorem 13]{hs}). 
More sufficient conditions guaranteeing stabilization at $E_2$ can be found in 
\cite[Theorem 3]{barnes}, \cite[Corollaries 1.5, 1.6 and Theorem 1.7]{dp1}, 
and \cite[Theorem 1.2]{dp2}.

The main result (Theorem 0.2) of \cite{z1} can be salvaged in the following way.
The above-mentioned general erroneous claim was used once, in the proof of 
Lemma 2.1. Lemma 2.1, in its turn, was used in the proof of Theorem 0.2 twice: 
first time (in the proof of Lemma 2.4) in the situation where $S$ is semisimple,
what is covered by the Hochschild--Serre result mentioned above,
and the second time in the situation where $L = S + I$ is the direct sum of 
algebras. In the latter case, one can invoke the K\"unneth theorem instead of 
Proposition 1.3.

I am grateful to Donald Barnes who pointed out this error.